\documentclass[a4paper,10pt]{article}

\usepackage{amsmath,enumerate,amsfonts,amssymb,amsthm}

\theoremstyle{plain}
\newtheorem*{Theo}{\bf Theorem}
\newtheorem{TheoN}{\bf Theorem}
\newtheorem*{Prop}{\bf Proposition}
\newtheorem{PropN}{\bf Proposition}

\newtheorem*{Ex}{\bf Example}

\theoremstyle{definition}
\newtheorem*{Def}{\bf Definition}

\theoremstyle{remark}
\newtheorem*{Rem}{\bf Remark}

 \DeclareMathOperator{\im}{Im}
\DeclareMathOperator{\Ker}{Ker} \DeclareMathOperator{\hess}{hess}
\DeclareMathOperator{\pf}{pf}

\DeclareMathOperator{\Hess}{Hess}
\begin{document}

\title{\bfseries Complex solutions \\ of Monge-Amp\`ere equations}
\author{Bertrand Banos}
\date{}
\maketitle

\footnotetext[1]{LMAM, Universit\'e de Bretagne Sud, Centre Yves Coppens, Campus de Tohannic\\
BP 573, 56017 VANNES, FRANCE\\ email:\texttt{bertrand.banos@univ-ubs.fr}}

\begin{abstract}
We describe a method to reduce partial differential equations of Monge-Amp\`ere type in 4 variables to complex partial differential equations in 2 variables. To illustrate this method, we construct explicit holomorphic solutions of the special lagrangian equation, the real Monge-Amp\`ere equations and the Plebanski equations.
\end{abstract}

\section*{Introduction}

A Monge-Amp\`ere equation is a partial differential equation which is non linear in a very specific way: its nonlinearity is the determinant one. In two variables, Monge-Amp\`ere equations are
$$
A\frac{\partial^2 f}{\partial q_1^2}+2B\frac{\partial^2
f}{\partial q_1\partial q_2}+C\frac{\partial^2 f}{\partial q_2^2}+
D\Big(\frac{\partial^2 f}{\partial q_1^2}\cdot\frac{\partial^2
f}{\partial q_2^2}-(\frac{\partial^2 f}{\partial q_1\partial
q_2})^2\Big)+E=0 \;\;\;\;\;,
$$
where coefficients $A$, $B$, $C$ and $D$ are smooth functions on jet space $J^1\mathbb{R}^2$. An important subfamily is the family of ``symplectic'' Monge-Amp\`ere equations, when coefficients are assumed to be smooth functions on cotangent space $T^*\mathbb{R}^2$. Contact geometry is then replaced by symplectic geometry.  In this work, we will only consider  Monge-Amp\`ere with constant coefficients which fall in this symplectic subfamily.

In higher dimensions, a Monge-Amp\`ere equation is a linear combination of the minors of the hessian matrix
$$
\Hess(f)=\Big(\frac{\partial^2 f}{\partial q_j\partial q_k}\Big)_{j,k=1\ldots n}\;\;\;\;\;.
$$
As examples in dimension $4$, we will cite the famous special lagrangian equation described  by Harvey and Lawson, or Plebanski equations and Grant equation obtained by reduction of Yang-Mills equations and Einstein equations:
$$
\begin{cases}
\displaystyle \frac{\partial^2 f}{\partial x_1\partial x_3}\frac{\partial^2
  f}{\partial x_2\partial x_4}- \frac{\partial^2 f}{\partial
  x_1\partial x_4}\frac{\partial^2 f}{\partial x_2\partial
  x_3}=1&\text{(Plebanski I equation)}\\
  &\\
\displaystyle \frac{\partial^2 f}{\partial x_1^2}\frac{\partial^2}{\partial x_2^2}-
\Big(\frac{\partial^2 f}{\partial x_1\partial x_2}\Big)^2+ \frac{\partial^2
  f}{\partial x_2\partial x_4}-\frac{\partial^2 f}{\partial
  x_1\partial x_3}=0&\text{(Plebanski II equation)}\\
 &\\
  \displaystyle \frac{\partial^2 f}{\partial x_1^2}+ \frac{\partial^2 f}{\partial
  x_1\partial x_2}\frac{\partial^2 f}{\partial x_3\partial x_4} -
\frac{\partial^2 f}{\partial x_1\partial x_4}\frac{\partial^2
  f}{\partial x_2\partial x_3}=0&\text{(Grant equation)}.\\
\end{cases}
$$

After the seminal paper of Lychagin (\cite{L}), geometry of Monge-Amp\`ere equations is quite well understood in 2 and 3 variables (\cite{LR},\cite{Ba1}, \cite{Ba2}, \cite{Ba3}) but dimension 4 remains mysterious. This is actually the first dimension in which the space of equivalent classes (modulo a change of independent and dependent variables)  is a real moduli space and discrete classification is not possible anymore. We propose here a method to reduce these 4 dimensional equations to complex equations in two variables.

In the first section, we recall the Lychagin correspondence between Monge-Amp\`ere equations and effective forms on the phase space and we describe briefly  classification results in dimension 2 and 3.
In the second section, we assume that the phase space is endowed with an extra compatible complex structure and we define the notion of  complex solution. We explain why complex solutions of a given Monge-Amp\`ere equation depend only of  its bieffective part.
In the third section we use this method to construct explicit complex solutions of the special lagrangian equation, the real Monge-Amp\`ere equations and the Plebanski equations.

\section{Monge-Amp\`ere operators and differential forms}

\subsection{An example: the special lagrangian equation}

Let $\displaystyle \Omega= \dfrac{i}{2} \big(dz_1\wedge d\overline{z_1}+\ldots + dz_n\wedge d\overline{z_n}\big)$ be the canonical K\"ahler form on  $\mathbb{C}^n$ and $\alpha=dz_1\wedge \ldots \wedge dz_n$ be the complex volume form. A special lagrangian submanifold is a real $n$-submanifold $L$ which is lagrangian with respect to $\Omega$ and which satisfies the special condition
$$
\im(\alpha)|_L=0.
$$

These submanifolds, introduced by Harvey and Lawson in their famous article \emph{Calibrated Geometries} (\cite{HL}) are minimal submanifolds of $\mathbb{C}^n$, and more generally of Calabi-Yau manifolds,  transverse in some sense to complex submanifolds. They have been extensively studied after the construction proposed by Strominger, Yau and Zaslow (\cite{SYZ}) of mirror partners of Calabi-Yau manifolds based upon an hypothetic special lagrangian fibration.

Some examples have been given by many people. We can cite for example Harvey and Lawson (\cite{HL}), Joyce (\cite{J1}, \cite{J2}, \cite{J3}, \cite{J4} \cite{J5}) and Bryant (\cite{Br1}, \cite{Br2}, \cite{Br3}).

For every smooth function $f$ on $\mathbb{R}^n$, the graph $$L_f=\Big\{(q+i\frac{\partial f}{\partial q}), \; q\in \mathbb{R}^n\Big\}$$
is a lagrangian submanifold of $\mathbb{C}^n$. The special lagrangian condition becomes then a differential equation on $f$, called the special lagrangian equation:
\begin{enumerate}[$\bullet$]
\item $n=2$: $\Delta f =0$
\item $n=3$: $\Delta f - \hess f =0$
\item $n=4$: $\Delta f - \hess_1 f - \hess_2 f -\hess_3 f - \hess_4 f=0$
\end{enumerate}
with $\Delta f$ the Laplace operator, $\hess f$ the determinant of the hessian matrix and $\hess_i f$ the $(i,i)$-minor of the hessian matrix.

The Monge-Amp\`ere operators theory developed by Lychagin (\cite{L}) generalizes this correspondence between ``calibrated'' lagrangian submanifolds of $\mathbb{R}^{2n}$ and Monge-Amp\`ere equations on $\mathbb{R}^n$.

\subsection{The Monge-Amp\`ere operators theory}

Let $M$ a $n$-dimensional manifold and $T^*M$ its cotangent bundle endowed with the symplectic canonical form $\Omega\in \Omega^2(T^*M)$. Denote by $q=(q_1,\ldots, q_n)$ a coordinates system on $M$ and $(q,p)$ the corresponding Darboux coordinates system on $T^*M$ such that
$$
\Omega = dq_1\wedge dp_1 + \ldots + dq_n\wedge dp_n.
$$

Let $\omega\in \Omega^n(T^*M)$ be a $n$ - differential form on the $2n$ dimensional manifold $T^*M$. The Monge-Amp\`ere operator $\Delta_\omega: C^\infty(M)\rightarrow \Omega^n(M)$ is defined by
$$
\Delta_\omega(f)=(df)^*(\omega)
$$
where $df: M\rightarrow T^*M$ is the differential of the smooth function $f$.

\begin{enumerate}[a)]

\item  A regular solution of the MA equation $\Delta_\omega=0$ is a smooth function $f$ on $M$ such that $\Delta_\omega(f)=0$.

\item A  generalized solution is a lagrangian submanifold $L$ of $(T^*M,\Omega)$ on which vanishes the form $\omega$:
$$\Omega|_L=0 \;\;\;\text{ and } \;\;\;  \omega|_L=0.$$

Note that a lagrangian submanifold of $T^*\mathbb{R}^n$ which projects isomorphically on $\mathbb{R}^n$ is the graph of an exact form $df:\mathbb{R}^n\rightarrow T^*\mathbb{R}^n$. Hence, a generalized solution can be thought of as a smooth patching of local regular solutions.

\item  Two MA equations $\Delta_{\omega_1}=0$ and $\Delta_{\omega_2}=0$ are said (locally) equivalent if there exists a (local) symplectomorphism $F: (T^*M,\Omega)\rightarrow (T^*M,\Omega)$ such that
$$
F^*(\omega_1)=\omega_2.
$$
The symplectomorphism $F$ transforms a generalized solution of $\Delta_{\omega_2}=0$ into a generalized solution of $\Delta_{\omega_1}=0$ but regular solutions are not preserved.
\end{enumerate}

For any $(n-2)$-form $\theta$, the equations $\Delta_{\omega}=0$ and $\Delta_{\omega+\theta\wedge\Omega}=0$ have the same solutions. We need also to introduce  effective $n$-forms, which satisfy
$$
\Omega\wedge\omega=0
$$
and the so-called Hodge-Lepage-Lychagin theorem (\cite{L}) establishes a one-to-one correspondence between effective forms and MA operators:

\begin{Theo}[Hodge-Lepage-Lychagin\\]

\begin{enumerate}[a)]
\item Every $n$-form $\omega$ can be uniquely decomposed
$$
\omega=\omega_0+\omega_1\wedge\Omega\;\;\;\; \text{with $\omega_0$ effective}
$$
\item Two effective $n$-forms which vanish on the same lagrangian subspaces are proportional.
\end{enumerate}
\end{Theo}

\begin{Ex}
Consider the $2$-dimensional equation $$\hess f =1.$$ The corresponding effective form is $\omega=dp_1\wedge dp_2 - dq_1\wedge dq_2$ which is transformed  into $\theta=dp_1\wedge dq_2+ dq_1\wedge dp_2$ by the partial Legendre transformation
$$
\Phi: (q_1,q_2,p_1,p_2)\mapsto (q_1,p_2,p_1,-q_2)
$$
Hence, $\hess f =1$ is equivalent to the Laplace equation $\Delta f= 0$. Choose then any harmonic function, for example $f(q_1,q_2)=e^{q_1}\cos(q_2)$. We obtain a generalized solution
$$L=\Big\{(q_1,-e^{q_1}\sin(q_2),e^{q_1}\cos(q_2),-q_2);\; (q_1,q_2)\in\mathbb{R}^2\Big\}.$$
which is, on an open subset,  the graph of the non trivial regular solution
$$
u(t_1,t_2)=t_2\arcsin(t_2e^{-t_1})+\sqrt{e^{2t_1}-t_2^2}.
$$
We know from J\"orgens theorem that this regular solution can not be defined on the whole plane.
\end{Ex}

\subsection{Classification results for $n=2$ and $n=3$}

In 1874, Sophus Lie raised the question of linearization of Monge-Amp\`ere equations: when a given equation is equivalent to a linear one ? In the formalism of Monge-Amp\`ere operators, this problem turns into a problem of the Geometric Invariant Theory. Studying action of the symplectic group on effective forms, a complete classification of equivalence classes of Monge-Amp\`ere equations (with constant coefficients) has been obtained  in \cite{LR} for $n=2$ and in \cite{LR} and \cite{Ba2} for $n=3$.

For $n=2$, every MAE with constant coefficient is linearizable (table \ref{table1}).

\begin{table}[!ht]
\begin{center}
\begin{tabular}{|c|c|c|c|}
\hline
\mathversion{bold} $\Delta_\omega=0$&\mathversion{bold} $\omega$ & \mathversion{bold} $\pf(\omega)$& \bf{Geometry}\\
\hline \hline
$\Delta f=0$ & $dq_1\wedge dp_2 - dq_2\wedge dp_1$&$1$& $A_\omega^2=-1$\\
\hline
$\square f=0$& $dq_1\wedge dp_2 + dq_2\wedge dp_1$&$-1$&$A_\omega^2=1$\\
\hline $\frac{\partial^2 f}{\partial q_1^2}=0$ & $dq_1\wedge dp_2$
& $0$&$A_\omega^2=0$\\
 \hline
\end{tabular}
\end{center}
\caption{Classification of 2-dimensional MAE}\label{table1}
\end{table}

Here, the pfaffian $\pf(\omega)$ is the scalar defined by $\omega\wedge\omega = \pf(\omega)\;\Omega\wedge\Omega$
and $A_\omega$ is the tensor defined by $\omega(\cdot\; ,\; \cdot )=\Omega(A_\omega\cdot\; ,\; \cdot )$. It is a complex structure for the elliptic equation and a product structure for the hyperbolic equation. It is explained in \cite{Ba3} how this unifying geometry co\"\i ncides, in the particular dimension $n=2$, with the famous generalized complex geometry introduced by Hitchin (\cite{Hi2}) and Gualtieri (\cite{Gu}).

For $n=3$, there are three non linear Monge-Amp\`ere equation: the real one, the special lagrangian one and the pseudo special lagrangian one (table \ref{table2}).

\begin{table}[!ht]
\begin{center}
\begin{tabular}{|c|c|c|c|}
\hline
&\mathversion{bold}$\Delta_\omega=0$ &
\mathversion{bold}$\varepsilon(g_\omega)$ & \mathversion{bold}$A_\omega$\\
\hline
\hline
 1& $\hess(f)=1$ &$(3,3)$&$A_\omega^2=1$\\
\hline
 2& $\Delta f- \hess(f)=0$ &$(0,6)$&$A_\omega^2=-1$\\
\hline
3& $\square f +\hess(f)=0$ &$(4,2)$&$A_\omega^2=-1$\\
\hline
\hline
 4& $\Delta f=0$ &$(0,3)$&$A_\omega^2=0$\\
\hline
5&$\square f=0$ &$(2,1)$&$A_\omega^2=0$\\
\hline
6&$\Delta_{q_2,q_3} f=0$&$(0,1)$&$A_\omega^2=0$\\
\hline
7&$\square_{q_2,q_3} f=0$&$(1,0)$&$A_\omega^2=0$\\
\hline
8&$\frac{\partial^2 f}{\partial q_1^2}=0$&$(0,0)$&$A_\omega^2=0$\\
\hline
\end{tabular}
\end{center}
\caption{Classification of 3-dimensional MAE}\label{table2}
\end{table}

In this table, $\varepsilon(g_\omega)$ is the signature of the Lychagin-Rubtsov metric (see \cite{LR}) defined by
$$
g_\omega(X,Y)\Omega^3= \iota_X(\omega)\wedge\iota_Y(\omega)\wedge\Omega,
$$
and $A_\omega$ is the Hitchin tensor (\cite{LR}), defined for effective $3$-forms by
$$
g_\omega(A_\omega\cdot\;,\;\cdot)=\Omega(\cdot\;,\;\cdot).
$$
It is explained in \cite{Ba1} how these invariants define a geometry of real or complex Calabi-Yau type.

For dimensions 2 and 3, the quotient space of Monge-Amp\`ere equations for  the action of symplectic linear group $G=Sp(2n,\mathbb{R})$ is thus a discrete space. More generally, a Monge-Amp\`ere equation corresponds to a conformal class of an effective form $\omega$. The orbit $G\cdot\omega$ of this form  is isomorphic to the quotient $G/G_\omega$ where $G_\omega=\{F\in G, F^*\omega=\omega\}$ is the stabilizer. For $n=4$,   the dimension of the group $G=Sp(8,\mathbb{R})$ is $36$, and dimension of the effective $4$-forms space is $42$. Dimensions of some stabilizers have been computed in \cite{DF} (table \ref{table3})

\begin{table}[!ht]
\begin{center}
\begin{tabular}{|c|c|}
\hline
{\bf Equation} & {\bf Stabilizer's dimension}\\
\hline
\hline
SLAG & 15\\
\hline
$\hess(f)=1$&15\\
\hline
Plebanski I & 13\\
\hline
Plebanski II & 14\\
\hline
Linear & $\geq$ 16\\
\hline
\end{tabular}
\end{center}
\caption{Stabilizers of some 4-dimensional MAE}\label{table3}
\end{table}

Moreover, generic effective forms have trivial stabilizer as it is explained in \cite{LR}. Hence, we need between $6$ and $21$ parameters to describe the quotient space around a Monge-Amp\`ere equation.

\section{Complex solutions and bieffective forms}

\subsection{Complex solutions}

A complex structure $\mathbb{J}$ on manifold $T^*M$ is said to be compatible with symplectic form $\Omega$, if $\Omega_{\mathbb{J}}= \Omega(\mathbb{J}\cdot\;,\;\cdot)$ is a $2$-form. The complex $2$-form $\Theta_\mathbb{J}=\Omega-i\Omega_{\mathbb{J}}$ is then a complex symplectic form.

In Darboux coordinates, such a compatible complex structure writes as
$$
\mathbb{J} = \begin{pmatrix}
A&B\\
C&A^t\\
\end{pmatrix} \;\;\;\;\text{ with } \begin{cases}
B^t=-B,\;\; C^t=-C&\\
A^2+BC=-1&\\
AB+BA^t=0&\\
AC+CA^t=0&\\
\end{cases}
$$
As in generalized complex geometry, there are two important families:
$$
\mathbb{J}=\begin{pmatrix}
A&0\\
0&A^t\\ \end{pmatrix} \;\;\; \text{ with $A$ a complex structure on $M$}
$$
and
$$
\mathbb{J}=\begin{pmatrix}
0&\theta\\
-\theta^{-1}&0\\ \end{pmatrix} \;\;\; \text{ with $\theta$ a symplectic form on $M$}
$$

\begin{Def}
Let $\Delta_\omega=0$ be a MAE on $M$ and let $\mathbb{J}$ be a compatible complex structure on $T^*M$. A $\mathbb{J}$-complex solution is a lagrangian and $\mathbb{J}$-complex submanifold $L$ of $T^*M$ on which $\omega$ vanishes:
$$
\Omega|_L=0\;; \;\;\;\;\;\;\; \mathbb{J}\;L=L\;; \;\;\;\;\;\;\; \omega|_L=0
$$
\end{Def}

It is worth mentioning that the condition ``complex lagrangian'' is equivalent to the condition ``real bilagrangian'' as it is proved in \cite{Hi3}.

\begin{Prop}[Hitchin] Let $(N,\Omega_1+i\Omega_2)$ be a complex symplectic manifold of complex dimension $2m$. A real $2m$-dimensional submanifold is a complex lagrangian submanifold if and only it is lagrangian with respect to $\Omega_1$ and $\Omega_2$.
\end{Prop}

\begin{Ex}
Let us identify $\mathbb{C}^4$ with $\mathbb{H}^2$ endowed with the three complex structures $I$, $J$ and $K$. We still denote by $\Omega$ the $I$ - K\"ahler form, $\Omega_J=\Omega(J\cdot\;,\;\cdot)$ and $\Omega_K=\Omega(K\cdot\;,\;\cdot)$. Then $\Omega-i\Omega_J$ is a $J$-complex symplectic form and it is well known that every $J$-complex lagrangian submanifold is special lagrangian.

In the formalism of Monge-Amp\`ere operators, this can be seen very simply. It is actually straightforward to check that the special lagrangian form writes as
$$
\im(\alpha)=\Omega_J\wedge\Omega_K
$$
and therefore, if $\Omega|_L=0$ and $\Omega_J|_L=0$ then $\im(\alpha)|_L=0$.
\end{Ex}

\subsection{Bieffective forms}

Two MAE $\Delta_{\omega}=0$ and $\Delta_{\omega+\theta_1\wedge\Omega+ \theta_2\wedge\Omega_{\mathbb{J}}}=0$ have the same $\mathbb{J}$-complex solutions. To understand this complex reduction, we need then to construct the bieffective part of $\omega$. This is the goal of this section.

Let $V$ be a complex symplectic space of real dimension $4m$ endowed with a complex symplectic form $\Theta=\Omega_1+i\Omega_2$. Denote by $\Lambda^k(V^*)$ the space of real $k$-forms on $V$ and $\Lambda^{p,q}(V^*)$ the space of $(p,q)$-complex forms, such that
$$
\Lambda^k(V^*)\otimes\mathbb{C}=\underset{p+q=k}{\bigoplus} \Lambda^{p,q}(V^*)
$$

Let us introduce for $j=1,2$ the operators $\top_j$ and $\bot_j$ defined by
$$
\begin{cases}
\top_j \theta = \theta\wedge\Omega_j&\\
&\\
\bot_j \theta = \iota_{X_{\Omega_j}}(\theta) & \text{ with $X_{\Omega_j}$ the unique bivector satisfying $\Omega_j(X_{\Omega_j})=1$}\\
\end{cases}
$$

They have the following properties (\cite{L})

\begin{enumerate}[a)]
\item $\bot_j: \Lambda^k(V^*)\rightarrow \Lambda^{k-2}(V^*)$ is into for $k\geq 2m+1$
\item $\top_j: \Lambda^k(V^*)\rightarrow \Lambda^{k+2}(V^*)$ is into for $k\leq 2m-1 $
\item $[\bot_j,\top_j] (\theta) = (2m-k)\theta$ for $\theta\in\Lambda^k(V^*)$.
\end{enumerate}

A $k$-form $\theta$ is said to be $\Omega_j$-effective if $\bot_j\theta=0$. For $k=2m$, this is equivalent to $\top_j\theta=0$.

Let $H=[\bot_1,\top_1]=[\bot_2,\top_2]$ and $M=[\bot_2,\top_1]$. We get then the complete list of so-called Lichnerowicz operators which satysfy the following (see \cite{Bo}):

\begin{Prop}[Verbitskii - Bonan]
$$
\begin{aligned}
\lbrack \bot_1,\top_1\rbrack &=H&&\lbrack \bot_2,\top_2\rbrack =H\\
\lbrack \bot_1,\top_2\rbrack &=-M&&\lbrack \bot_2,\top_1\rbrack =M\\
\lbrack \bot_1,\bot_2\rbrack &=0&&\lbrack \top_1,\top_2\rbrack =0\\
\lbrack \bot_1,H\rbrack &=-2\bot_1&&\lbrack \bot_2,H\rbrack =-2\bot_2\\
\lbrack \top_1,H\rbrack &=2\top_1&&\lbrack \top_2,H\rbrack =2\top_2\\
\lbrack \bot_1,M\rbrack &=-2\bot_2&&\lbrack \bot_2,M\rbrack =2\bot_1\\
\lbrack \top_1,M\rbrack &=-2\top_2&&\lbrack \top_2,M\rbrack =2\top_1\\
&&\;\;\lbrack H,M\rbrack =0&\\
\end{aligned}
$$
\end{Prop}

We obtain then a representation of the Lie algebra $sl(2,\mathbb{C})\otimes sl(2,\mathbb{C})$ on $\Lambda^*(\mathbb{V^*})\otimes\mathbb{C}$, defining
$$
\begin{cases}
E_1=\frac{1}{2}(\bot_1+i\bot_2)&\\
F_1=\frac{1}{2}(\top_1-i\top_2)&\\
H_1=\frac{1}{2}(H+iM)&
\end{cases}\;\;\;\; \begin{cases} E_2=\frac{1}{2}(\bot_1-i\bot_2)&\\F_2=\frac{1}{2}(\top_1+i\top_2)&\\H_2=\frac{1}{2}(H-iM)\\ \end{cases}
$$

Representation theory of Lie algebras gives us the existence and uniqueness of bieffective part of a $2m$-form.

\begin{TheoN}[\cite{Ba4}]
Every $2m$-form $\omega\in \Lambda^{2m}(V^*)$ can be decomposed into a sum
$$
\omega=\omega_0+ \omega_1\wedge\Omega_1+\omega_2\wedge\Omega_2
$$
with $\omega_0$ bieffective, that is $\omega_0\wedge \Omega_1=0=\omega_0\wedge\Omega_2$.
Moreover the bieffective part $\omega_0$ is unique.
\end{TheoN}

\begin{proof}
From Weyl's theorem, we know that $\Lambda^*(\mathbb{V^*})\otimes\mathbb{C}$ decomposes as a unique direct sum of irreducible subspaces. Let $W$ such an irreducible subspace. Since $H_1$ and $H_2$ commute, they admit a common eigenvector $x\in W$. But $E_1^pE_2^qx$ is also a common eigenvector. Their exist then $p$ and $q$ such that $z=E_1^pE_2^qx$ is a primitive vector, that is $E_1z=E_2z=0$. Therefore, $W=\mathcal{G}z$ and every vector $w$ in $W$ writes as
$$
w=\sum_{j,k} a_{jk}F_1^jF_2^k z
$$
We deduce that every $2m$-form can be uniquely decomposed into a finite sum
$$
\omega=\sum_{j,k} (\alpha_{jk}+i\beta_{jk})\wedge(\Omega_1+i\Omega_2)^j\wedge(\Omega_1-i\Omega_2)^k
$$
with $\alpha_{jk}$ and $\beta_{jk}$ primitive. Noting now that $2m$- primitive forms are the bieffective forms we obtain the result.
\end{proof}

We give now an explicit formula for $4m=8$, which can be easily implemented on a computer.

\begin{Prop}
In dimension $8$, the bieffective par $\omega_0$ of a $4$-form $\omega$ is
$$
\omega_0=\theta-\frac{1}{4}\Big\{\top_2\bot_2\theta+\top_1\bot_1\theta
-\frac{1}{4}M(M\theta-\top_1\bot_2\theta+\top_2\bot_1\theta)\Big\}
$$

where

$$
\theta=\omega-\frac{(3\bot_1^2\omega-\bot_2^2\omega)}{64}\Omega_1^2
-\frac{\bot_1\bot_2\omega}{8}
-\frac{(3\bot_2^2\omega-\bot_1^2\omega)}{64}\Omega_2^2
$$
\end{Prop}

\begin{proof}
We know that
$$
\omega=\omega_0+\omega_1\wedge\Omega_1+\omega_2\wedge\Omega_2+
\omega_{11}\Omega_1\wedge \Omega_1+
\omega_{12}\Omega_1\wedge\Omega_2+\omega_{22}\Omega_2\wedge\Omega_2,
$$
with $\omega_0$, $\omega_1$ and $\omega_2$ primitive. Using Verbitskii-Bonan relations, we obtain
$$
\bot_1\omega=2\omega_1-M\omega_2+(6\omega_{11}+2\omega_{22})\Omega_1 +
2\omega_{12}\Omega_1
$$
and then
$$
\bot_1^2\omega=24\omega_{11}+8\omega_{22} \;\;\;  \bot_2\bot_1\omega=8\omega_{12}.
$$
Starting from $\bot_2\omega$ we obtain also
$$
\bot_2^2\omega=8\omega_{11}+24\omega_{22}.
$$
Therefore,
$$
\omega_{11}=\frac{3\bot_1^2\omega-\bot_2^2\omega}{64},\;\;\;\omega_{22}=\frac{3\bot_2^2\omega-\bot_1^2\omega}{64},\;\;\; \omega_{12}=\frac{\bot_1\bot_2\omega}{8}.
$$

Define now $\theta=\omega_0+\omega_1\wedge\Omega_1+\omega_2\wedge\omega_2$. Since $\bot_1\theta=2\omega_1-M\omega_2$ and $\bot_2\theta=M\omega_1+2\omega_2$, we deduce that
$$
M\top_1\omega_1=[M,\top_1]\omega_1+\top_1M\omega_1=2\top_2\omega_1+\top_1(\bot_2\theta-2\omega_2)=2\top_2\omega_1-2\top_1\omega_2+\top_1\bot_2\theta
$$
and similarly
$$
M\top_2\omega_2=2\top_2\omega_1-2\top_1\omega_2-\top_2\bot_1\theta.
$$
and therefore
$$
M\theta=M\top_1\omega_1+M\top_2\omega_2=4(\top_2\omega_1-\top_1\omega_2)+\top_1\bot_2\theta -
\top_2\bot_1\theta.
$$
Moreover, the computation of $M(\top_2\omega_1-\top_1\omega_2)$ gives
$$
M(\top_2\omega_1-\top_1\omega_2)=-4(\top_1\omega_1+\top_2\omega_2)+\top_2\bot_2\theta+\top_1\bot_1\theta.
$$
Finally,
$$
4(\omega_1\wedge\Omega_1+\omega_2\wedge\Omega_2)=
\top_2\bot_2\theta+\top_1\bot_1\theta
-\frac{M}{4}(M\theta-\top_1\bot_2\theta + \top_2\bot_1\theta).
$$
\end{proof}

\subsection{Action of the complex symplectic group}

For simplicity, we restrict now to complex dimension $4$: $V$ is a $4$-dimensional complex vector space endowed with a complex symplectic form $\Theta=\Omega_1+i\Omega_2$. The space of real bieffective $4$-forms is
$$
\Lambda_{BE}^4(V^*)=\big\{\omega\in \Lambda^4(V^*), \omega \wedge\Omega_1=\omega\wedge\Omega_2=0\big\}.
$$

Let $\Lambda^{2,0}_0(V^*)$ be the $5$ - dimensional complex vector space of $(2,0)$-complex forms which are effective with respect to $\Theta=\Omega_1+i\Omega_2$:
$$
\Lambda^{2,0}_0(V^*)=\big\{\theta\in \Lambda^{2,0}(V^*),\; \theta\wedge\Theta=0\big\}
$$
The exterior product is non degenerate on $\Lambda^{2,0}_0(V^*)$ and it defines a non degenerate symmetric inner product
$$
<\theta_1,\theta_2>\;\Theta^2=\theta_1\wedge\theta_2\;\;.
$$

We are going to identify bieffective forms with hermitian forms on $\Lambda^{2,0}_0(V^*)$.

\begin{Def}
The hermitian form $Q_\omega$ on $\mathbb{C}^5=\Lambda^{2,0}_0(V^*)$ associated with a bieffective $4$-form $\omega$ is:
$$
Q_\omega(\theta_1,\theta_2)\; (\Theta\wedge\bar{\Theta})^2 = \omega\wedge \theta_1 \wedge \overline{\theta_2}\;\; .
$$
\end{Def}

\begin{TheoN}
The map
$$\begin{aligned} \Lambda^4_{BE}(V^*)&\rightarrow su(5)\\ \omega &\mapsto Q_\omega\end{aligned}$$ is an isomorphism.

Moreover, the group $Sp(4,\mathbb{C})/\mathbb{Z}_2$ identifies with $SO(5,\mathbb{C})$ and its action on $\Lambda^4_{BE}(V^*)$ is the Hermite action of $SO(5,\mathbb{C})$ on $su(5)$.
\end{TheoN}

\begin{proof}

We prove first that $\displaystyle \Lambda_{BE}^4(V^*)=\Lambda^{2,0}_0(V^*)\otimes \overline{\Lambda^{2,0}_0(V^*)}$.
Let $\omega\in \Lambda_{BE}^4(V^*)$ and consider its decomposition
$$
\omega=\omega_{40}+\omega_{31}+\omega_{22}+\omega_{13}+\omega_{04}
$$
with $\omega_{pq}\in \Lambda^{p,q}(V^*)$. Since $\bot_{\mathbb{C}}: \Lambda^{p,q}(V^*)\rightarrow \Lambda^{p-2,q}(V^*)$ is into for $p\geq 3$ and $\overline{\bot_{\mathbb{C}}}: \Lambda^{p,q}(V^*)\rightarrow \Lambda^{p,q-2}(V^*)$ is into for $q\geq 3$, and since $\bot_{\mathbb{C}}\omega=\overline{\bot_{\mathbb{C}}}\omega=0$, we deduce that $\omega\in \Lambda^{2,2}$. But
$$
\Ker\Big(\bot_{\mathbb{C}}:\Lambda^{2,2}\rightarrow \Lambda^{0,2}\Big) = \Lambda^{2,0}_0(V^*)\otimes \Lambda^{0,2}(V^*)
$$
so
$$
\Ker\Big(\bot_{\mathbb{C}}\Big) \cap \Ker \Big(\overline{\bot_{\mathbb{C}}}\Big) = \Lambda^{2,0}_0(V^*)\otimes \overline{\Lambda^{2,0}_0(V^*)}.
$$

We deduce that $\displaystyle \Lambda_{BE}^4(V^*)$ and $su(5)$ have same dimension. Since $\omega\mapsto Q_\omega$ is injective, this is an isomorphism.

Now, the action of $Sp(4,\mathbb{C})$ preserves this symmetric product, with kernel $\mathbb{Z}_2$. Since $dim_{\mathbb{C}}(Sp(4,\mathbb{C}))=10=dim_{\mathbb{C}}(SO(5,\mathbb{C}))$, we deduce that
$$
Sp(4,\mathbb{C})/\mathbb{Z}_2 = SO(5,\mathbb{C}).
$$
Moreover, we have
$$
Q_{F^*\omega}(\theta_1,\theta_2)(\Theta\wedge\bar{\Theta})^2 = F^*(\omega)\wedge \theta_1 \wedge \overline{\theta_2}=\omega\wedge (F^{-1})^*(\theta_1)\wedge(F^{-1})^*(\overline{\theta_2})
$$
so
$$
Q_{F^*\omega}=\overline{F^{-1}}^t Q_\omega F^{-1}.
$$
\end{proof}

This action is completely described by Hong in \cite{Ho}. Let us briefly explain this result. Let $Q$ be a hermitian matrix and define the canonical form of $Q$ as a direct sum of three hermitian matrices:
$$
J(Q)= H_P(Q)\oplus K_N(Q)\oplus K_C(Q)\;\;,
$$
which are obtained from Jordan blocks of $Q$ as follows:
\begin{enumerate}[a)]
\item $\displaystyle H_P(Q)= H_{m_1}(\lambda_1)\oplus \ldots \oplus H_{m_p}(\lambda_q)$, where all $\lambda_i\geq 0$ and $\lambda_i^2$ are the positive eigenvalues of $QQ^t$.
\item $\displaystyle K_N(Q)= K_{2n_1}(\mu_1)\oplus \ldots \oplus H_{2n_r}(\mu_r)$, where all $\mu_i>0$ and $-\mu_i^2$ are the negative eigenvalues of $QQ^t$.
\item $\displaystyle K_C(Q)= L_{2k_1}(\xi_1)\oplus \ldots \oplus L_{2k_s}(\xi_s)$, where $\xi_i^2$ are the non real  eigenvalues of $QQ^t$,
\end{enumerate}
with
$$
K_{2n}(\mu)=\begin{pmatrix} 0 & -iH_n(\mu)\\iH_n(\mu)& 0\end{pmatrix},\;\;\;\;
L_{2k}(\xi)=\begin{pmatrix} 0 & H_k(\xi)\\H_k^\star(\xi)& 0\end{pmatrix}
$$
and
$$
2H_m(\lambda)=\begin{pmatrix} 0 & 0 & \ldots & 0 & 1 & 2\lambda\\0 &\ldots & 0 & 1 & 2\lambda & 1\\\vdots&\ldots&1&2\lambda&1&0\\
\vdots          &1&2\lambda& 1&\ldots&\vdots\\ 1&2\lambda&1&\ldots & 0 & 0\\2\lambda&1&0&\ldots & 0 & 0\\
\end{pmatrix} + i\begin{pmatrix} 0&1&0&\ldots & 0 & 0\\ -1&0&1&\ldots & 0 & 0\\\vdots&-1&0& 1&\ldots&\vdots\\\vdots&\ldots&-1&0&1&0\\
0 &\ldots & 0 & -1 & 0 & 1\\ 0 & 0 & \ldots & 0 & -1 & 0\\\end{pmatrix}
$$
\begin{Theo}[Hong]
Let $Q$ be a hermitian matrix. Then there exists $F$ complex orthogonal and $\varepsilon=(\varepsilon_1,\ldots,\varepsilon_p)$ with $\varepsilon_i=\pm 1$ such that $\overline{F}^t Q F= J^\varepsilon(Q)$ with
$$
J^\varepsilon(Q)= H^\varepsilon_P(Q)\oplus K_N(Q)\oplus K_C(Q)=  \Big(\varepsilon_1H_{m_1}(\lambda_1)\oplus \ldots \oplus \varepsilon_p H_{m_p}(\lambda_q)\Big)\oplus K_N(Q)\oplus K_C(Q)
$$
\end{Theo}

It is therefore difficult to give a complete classification of all possible complex reductions. Nevertheless, to characterize the orbit of a bieffective form $\omega$, it is necessary to know
\begin{enumerate}[a)]
\item the signature $\varepsilon(Q_\omega)$ of $Q_\omega$,
\item the spectrum of $Q_\omega Q_\omega^t$.
\end{enumerate}
and it will be sufficient for the examples we are interested in.

\subsection{The complex lagrangian grassmannian}

Denote by $\text{Gr}_\omega$ the set of all complex lagrangian planes of the complex symplectic space $(V,\Theta)$ on which vanishes the bieffective form $\omega$.

The complex isomorphism $\Theta: V\rightarrow \Lambda^{1,0}(V^*)$ transforms a complex basis of such a plane $L$ into a decomposable effective $(2,0)$ forms $\theta_L\in \Lambda^{2,0}_0(V^*)$. Note that $\theta_L$ is decomposable if and only if $\theta_L\wedge\theta_L=0$.

Moreover, the condition $\omega|_L=0$ is equivalent to the condition $\omega\wedge\theta_L\wedge\overline{\theta_L}=0$.

\begin{Prop}
The grassmannian $\text{Gr}_\omega$ identifies with the real algebraic subvariety of $P^4(\mathbb{C})=P(\Lambda^{2,0}_0)$ defined by
$$
\theta\wedge\theta=0\;\; \text{ and } \;\;\; Q_\omega(\theta)=0.
$$
\end{Prop}

\begin{Rem}
This grassmannian could be empty, for example if $Q_\omega$ is positive-definite. It could happen therefore that the Monge-Amp\`ere equation $\Delta_\omega=0$ has no generalized complex solution.
\end{Rem}

\section{Some examples in dimension 4}

\subsection{The choice of the complex structure}

The crucial point in this method is the choice of the compatible complex structure, which should depend on the initial Monge-Amp\`ere equation. We choose here five simple complex structures, and give a corresponding complex Darboux coordinates system $(z_1,z_2,u_1,u_2)$ in which
$$\Theta_\mathbb{J}=\Omega-i\Omega_\mathbb{J}=dz_1\wedge du_1+dz_2\wedge du_2$$
The initial coordinate systems on $T^*\mathbb{R}^4$ is still $(q,p)$ with
$$
\Omega=dq_1\wedge dp_1 +dq_2\wedge dp_2 + dq_3\wedge dp_3 + dq_4\wedge dp_4 .
$$
We define
$$
A=\begin{pmatrix} 0&-1& 0 & 0 \\ 1&0&0&0\\0&0&0&-1\\0&0&1&0\\ \end{pmatrix};\;\;\; \tilde{A}=\begin{pmatrix} 0&-1& 0 & 0 \\ 1&0&0&0\\0&0&0&1\\0&0&-1&0\\ \end{pmatrix} ;\;\;\; A_2=\begin{pmatrix} 1&-2& 0 & 0 \\ 1&-1&0&0\\0&0&1&-2\\0&0&1&-1\\\end{pmatrix}
$$
and
$$
J=\begin{pmatrix} A& 0\\ 0 & A^t \end{pmatrix} \;\;\;\;\;\;\; \begin{cases} z_1=q_1+iq_2 & \;\; u_1=p_1-ip_2\\ z_2=q_3+iq_4 & \;\; u_2= p_3-ip_4\end{cases}
$$
$$
K=\begin{pmatrix} 0& A\\ A & 0 \end{pmatrix} \;\;\;\;\;\;\; \begin{cases} z_1=q_1+ip_2 & \;\; u_1=iq_2+p_1\\ z_2=q_3+ip_4 & \;\; u_2= iq_4+p_3\end{cases}
$$
$$
\tilde{J}=\begin{pmatrix} \tilde{A}& 0\\ 0 & \tilde{A}^t \end{pmatrix} \;\;\;\;\;\;\; \begin{cases} z_1=q_1+iq_2 & \;\; u_1=p_1-ip_2\\ z_2=q_3-iq_4 & \;\; u_2= p_3+ip_4\end{cases}
$$
$$
\tilde{K}=\begin{pmatrix} 0& \tilde{A}\\ \tilde{A} & 0 \end{pmatrix} \;\;\;\;\;\;\; \begin{cases} z_1=q_1+ip_2 & \;\; u_1=iq_2+p_1\\ z_2=q_3-ip_4 & \;\; u_2= -iq_4+p_3\end{cases}
$$

$$
J_2=\begin{pmatrix} A_2& 0\\ 0 & A_2^t \end{pmatrix} \;\;\;\;\;\;\; \begin{cases} z_1=q_1+(-1+i)q_2 & \;\; u_1=(1-i)p_1-ip_2\\ z_2=q_3+(-1+i)q_4 & \;\; u_2= (1-i)p_3-ip_4\end{cases}
$$

\subsection{Examples}

We study now the special lagrangian equation, the two real Monge Amp\`ere equations $\hess f =\pm 1$, the two Plebanski equations and the Grant equation. The corresponding effective forms on $(T^*\mathbb{R}^4, \Omega)$ are:
$$
\begin{aligned}
\omega_{SLAG}& = \im\big (dq_1 + i dp_1 )\wedge (dq_2 + i dp_2)\wedge(dq_3 + i dp_3)\wedge (dq_4 + i dp_4)\big)\\
\omega_{H+}& =dp_1\wedge dp_2\wedge dp_3\wedge dp_4 - dq_1\wedge dq_2\wedge dq_3\wedge dq_4\\
\omega_{H-}& =dp_1\wedge dp_2\wedge dp_3\wedge dp_4 + dq_1\wedge dq_2\wedge dq_3\wedge dq_4\\
\omega_{PI}&= dq_1\wedge dq_2\wedge dp_1\wedge dp_2 - dq_1\wedge dq_2\wedge dq_3\wedge dq_4\\
\omega_{PII}&=dq_1\wedge dq_2+ dq_3\wedge dp_2 + dq_1\wedge dq_2\wedge dq_4\wedge dp_1+dq_3\wedge dq_4\wedge dp_1\wedge dp_2\\
\omega_G&=dq_2\wedge dq_3\wedge dq_4\wedge dp_1 - dq_1\wedge dq_3\wedge dp_1\wedge dp_3\\
\end{aligned}
$$

We compute  for the five compatible complex structures defined above the bieffective part  of these forms and the signature of  the hermitian form associated with. These invariants are given in  table \ref{table5} and have to be compared with invariants for simple complex equations given in table \ref{table4}.

In table \ref{table4}, $\phi$ is a holomorphic function in $(z_1,z_2)$ and we note $\displaystyle \phi_{jk}=\dfrac{\partial^2\phi}{\partial z_j\partial z_k}$.

\begin{table}[!ht]
\begin{center}
\begin{tabular}{|c|c|c|}
\hline
$\varepsilon(Q_\omega)$& $\Delta_\omega=0$ & spectrum$(Q_\omega Q_\omega^t)$\\
\hline
 & $|\phi_{11}|^2=0$ & (0,0,0,0,0)\\
(1,0)&&\\
 & $|\phi_{12}|^2=0$ & (1,0,0,0,0)\\
\hline
\hline
 & $|\phi_{11}|^2=1$ & (0,0,0,0,0)\\
&&\\
(1,1) & $|\phi_{12}|^2=1$ & (1,0,0,0,0)\\
&&\\
 & $|\phi_{11}|^2-|\phi_{22}|^2=0$ & (-1,-1,0,0,0)\\
\hline
 & $|\phi_{11}|^2+|\phi_{12}|^2=0$ & (1,0,0,0,0)\\
(2,0)&&\\
& $|\phi_{11}|^2+|\phi_{22}|^2=0$ & (1,1,0,0,0)\\
\hline
 & $|\phi_{11}|^2+|\phi_{12}|^2=1$ & (1,0,0,0,0)\\
(2,1)&&\\
 & $\phi_{12}+\overline{\phi_{12}}-|\phi_{11}|^2=0$ & (0,0,0,0,0)\\
\hline

\end{tabular}
\end{center}
\caption{Simple complex Monge-Amp\`ere equations in complex dimension 2}\label{table4}
\end{table}

\begin{table}[!ht]
\begin{center}
\begin{tabular}{|c||c|c|c|c|c|}
\hline
{\bf Equation} & $\bold{J}$ & $\bold{K}$ & $\bold{\tilde{J}}$ & $\bold{\tilde{K}}$ & $\bold{J_2}$\\
\hline
\hline
SLAG& 0&0&0&0&(1,1)\\
\hline
$\hess(f)=1$& (1,1) & 0 & (1,1) & 0 & (1,1)\\
\hline
$\hess(f)=-1$& (2,0) & (3,2) & (2,0) & (3,2)& (2,0)\\
\hline
Plebanski I & (2,0) & (3,2) & (1,1) & (3,2) & (2,0)\\
\hline
Plebanski II & (2,1) & (3,2) & (1,0) & (3,2) & (2,1)\\
\hline
Grant & (3,2)& (3,2) & (3,2) & (3,2) & (3,2)\\
\hline
\end{tabular}
\end{center}
\caption{Invariants for some Monge-Amp\`ere equations in dimension 4}\label{table5}
\end{table}

We see in table \ref{table5} that this method fails for the Grant equation, at least for this choice of complex structures: signature is always (3,2) and the corresponding complex equation is therefore fully non degenerate.  We study  the other cases in more details.

\subsubsection{Special lagrangian equation}

As, we have seen the special lagrangian form $\omega_{SLAG}$ has zero bieffective part for $J$ and for $K$ since
$$
\omega_{SLAG}=\Omega_{J}\wedge\Omega_{K}
$$
This is the best situation, since every complex lagrangian submanifold is solution. We note that this is the same for  $\tilde{J}$ and $\tilde{K}$ and more generally for any complex structure $\mathbb{J}= F^{-1}JF$ with $F$ in $SU(4)$, since
$$
\omega_{SLAG}=F^*\omega_{SLAG} = F^*\Omega_{J}\wedge F^*\Omega_K = \Omega_{F^{-1}JF}\wedge \Omega_{F^{-1}KF}
$$
Nevertheless, it does not give new solutions: we already know that the action of $F\in SU(n)$ transforms a special lagrangian submanifold into an other special lagrangian submanifold.

This is the reason of the choice of $J_2$: we were looking for a simple complex structure which is not in $so(8)$ in order to construct other examples of solutions.

In the complex Darboux coordinates system, the bieffective part is
$$
\begin{aligned}\omega_{SLAG}^{BE}=\frac{1}{8}\Big\{ &(1+2i)dz_1\wedge dz_2\wedge d\overline{z_1}\wedge d\overline{u_{2}} + (-1-2i)dz_1\wedge dz_2\wedge d\overline{z_2}\wedge d\overline{u_{1}}\\
&+(1-2i)dz_1\wedge du_2\wedge d\overline{z_1}\wedge d\overline{z_{2}} +(1+2i)dz_1\wedge du_2\wedge d\overline{u_1}\wedge d\overline{u_{2}}\\
&+(-1+2i)dz_1\wedge du_1\wedge d\overline{z_1}\wedge d\overline{z_{2}} + (-1-2i)dz_2\wedge du_1\wedge d\overline{u_1}\wedge d\overline{u_{2}}\\
&+(1-2i)du_1\wedge du_2\wedge d\overline{z_1}\wedge d\overline{u_{2}}+ (-1+2i)du_1\wedge du_2\wedge d\overline{z_2}\wedge d\overline{u_{1}}\Big\}
\end{aligned}
$$
and this is straightforward to check that
$$
\omega_{SLAG}^{BE}=\frac{\sqrt{5}}{4} \Big\{ dZ_1\wedge dU_2\wedge d\overline{Z_1}\wedge d\overline{U_{2}} - dZ_2\wedge dU_1\wedge d\overline{Z_2}\wedge d\overline{U_{1}}\Big\}
$$
where $(Z_1,Z_2,U_1,U_2)$ is the complex Darboux coordinates system
$$
\begin{aligned}
&Z_1=\frac{\alpha z_1+\alpha^{-1} u_1}{i\sqrt{2}} \;\; ; \;\; U_1=\frac{\alpha z_1-\alpha^{-1} u_1}{i\sqrt{2}}\\
&Z_2=\frac{\alpha z_2-\alpha^{-1} u_2}{\sqrt{2}} \;\; ; \;\; U_2=\frac{\alpha z_2+\alpha^{-1} u_2}{\sqrt{2}}\\
\end{aligned}
$$
$$
\alpha^2=\frac{1+2i}{\sqrt{5}}
$$

We obtain then the following result:

\begin{PropN}
Let $\phi$ be a holomorphic solution of
$$
\Big|\frac{\partial^2 \phi}{\partial z_1^2}\Big|^2 - \Big|\frac{\partial^2 \phi}{\partial z_2^2}\Big|^2=0
$$
and let $L_\phi$ be the submanifold $\displaystyle L_\phi=\big\{(z_1,z_2,\frac{\partial \phi}{\partial z_1},\frac{\partial \phi}{\partial z_2})\big\}\subset \mathbb{C}^4$.

Then $F^{-1}(L_\phi)$ is special lagrangian in $(T^*\mathbb{R}^4,\Omega,\omega_{SLAG})$ where $F(q,p)=(z,u)$ with
$$
\begin{aligned}
& z_1=\frac{\alpha q_1+(-1+i)\alpha q_2+ (1-i)\alpha^{-1} p_1 - i\alpha^{-1} p_2}{i\sqrt{2}}\\
&z_2=\frac{\alpha q_3 + (-1+i)\alpha q_4 - (1-i)\alpha^{-1} p_3+ i\alpha^{-1} p_4}{\sqrt{2}}\\
&u_1=\frac{\alpha q_1+(-1+i)\alpha q_2- (1-i)\alpha^{-1} p_1 + i\alpha^{-1} p_2}{i\sqrt{2}}\\
&u_2=\frac{\alpha q_3 + (-1+i)\alpha q_4 + (1-i)\alpha^{-1} p_3 - i\alpha^{-1} p_4}{\sqrt{2}}\\
\end{aligned}
$$
\end{PropN}

\subsubsection{Real Monge-Amp\`ere equations}

For every holomorphic function $\phi=f+ig: \mathbb{C}^2\rightarrow \mathbb{C}$, we have
$$
\hess_\mathbb{R} f = |\hess_\mathbb{C} \phi|^2
$$
so for $J$, $\tilde{J}$ and $J_2$, which come from complex structures on $\mathbb{R}^4$, the complex reduction of $\hess f = \pm 1$ are
$$
|\hess \phi|^2=\pm 1
$$
which are equivalent to
$$
\Big|\frac{\partial^2 \psi}{\partial z_1^2}\Big|^2 =\pm \Big|\frac{\partial^2 \psi}{\partial z_2^2}\Big|^2
$$

For example, for $J$, the corresponding symplectomorphism is the partial Legendre transform
$$
G(z_1,z_2,u_1,u_2)=(u_1,z_2,-z_1,u_2)
$$
We obtain then the  following results:

\begin{PropN}
Let $\phi$ be a holomorphic solution of
$$
\Big|\frac{\partial^2 \phi}{\partial z_1^2}\Big|^2 - \Big|\frac{\partial^2 \phi}{\partial z_2^2}\Big|^2=0
$$

Then $G^{-1}(L_\phi)$ is a generalized solution of $\hess f =1$ with
$$
G(q,p)=(p_1-ip_2,q_3+iq_4,-q_1-iq_2,p_3-ip_4)
$$
\end{PropN}

\begin{PropN}
Let $\phi$ a holomorphic function of the form $\phi(z_1,z_2)=a(z_1)b(z_2)$.
Then $G^{-1}(L_\phi)$ is a generalized solution of $\hess f =-1$.
\end{PropN}

Moreover, we see in table 5 that $\omega_{H+}$ has no bieffective part for $K$ and $\tilde{K}$. We get immediately Proposition 4:

\begin{PropN}
Any complex lagrangian surface in $(T^*\mathbb{R}^4,\Omega,K)$ or  $(T^*\mathbb{R}^4,\Omega,\tilde{K})$ is a generalized solution of $\hess f =1$.
\end{PropN}

This result is underlying the strong relationship there is between special lagrangian geometry and ``real special lagrangian'' geometry or ``split special lagrangian geometry'', as described in \cite{HL2}.

\subsubsection{Plebanski equations}

Computing the bieffective part of $\omega_{PI}$, we obtain that the $J$-complex reduction of Plebanski I equation is
$$
|\phi_{11}|^2=-1
$$
which is equivalent up the partial Legendre transform $G$  to
$$
|\psi_{11}|^2+|\psi_{12}|^2=0
$$

\begin{PropN}
Let $\phi$ a holomorphic function of the form $\phi(z_1,z_2)=az_1+b(z_2)$.
Then $G^{-1}(L_\phi)$ is a generalized solution of Plebanski I equation.
\end{PropN}

Similarly, the $\tilde{J}$-complex reduction is
$$
|\phi_{12}|^2=1
$$

\begin{PropN}
Let $\phi$ a holomorphic function on $(\mathbb{R}^4,\tilde{A})$ of the form $$\phi(z_1,z_2)=z_1z_2+a(z_1)+b(z_2).$$
Then its real part is a regular solution of Plebanski I equation.
\end{PropN}

Finally, same computations give analog results for Plebanski II equation:

\begin{PropN}
Let $\phi$ a holomorphic function on $(\mathbb{R}^4,A)$ of
$$
\phi_{12}+\overline{\phi_{12}}+|\phi_{11}|^2=0
$$
Then its real part is a regular solution of Plebanski II equation.
\end{PropN}

\begin{PropN}
Let $\phi$ a holomorphic function on $(\mathbb{R}^4,\tilde{A})$ of the form
$$\phi(z_1,z_2)=a(z_2)+b(z_2)z_1.$$
Then its real part is a regular solution of Plebanski II equation.
\end{PropN}

\section*{Conclusion}

Studying geometry of $4$-bieffective forms on $\mathbb{R}^8$, we have reduced important equations in Physic to simple - but non empty - complex equations.

This has been done for arbitraries complex structures.  It would be interesting now, for a given equation in four variables, to study all possible compatible complex structures and to parameterize in this way analytical solutions by pair of compatible complex structures on $\mathbb{R}^8$ , and holomorphic functions on $\mathbb{C}^2$. Note that bieffective part still exists for $2m$-forms on $\mathbb{R}^{4m}$ and such a parametrization should also exist for equations with $2m$ variables.

A more global approach would be also interesting. We have considered only complex structures with constant coefficients but our decomposition theorem remains valid on a complex symplectic manifold. This suggests to understand  ``Monge-Amp\`ere calibrations'' on complex symplectic manifolds as a generalization of the special lagrangian calibration on HyperK\"ahler manifolds.

\section*{Acknowledgements}

I am very grateful to Volodya Rubtsov and Jenya Ferapontov for  their helpful comments and suggestions.

\end{document}